\documentclass{amsart}

\usepackage{amssymb}
\usepackage{mathrsfs}

\theoremstyle{plain}
\newtheorem{theorem}{Theorem}[section]

\newtheorem{conjecture}[theorem]{Conjecture}
\newtheorem{proposal}[theorem]{Proposal}

\theoremstyle{definition}

\theoremstyle{remark}

\numberwithin{equation}{section}

\begin{document}
\newcommand{\Z}{\mathbb{Z}}
\newcommand{\Q}{\mathbb{Q}}
\newcommand{\R}{\mathbb{R}}
\newcommand{\C}{\mathbb{C}}
\newcommand{\divides}{\mid}
\newcommand{\doesnotdivide}{\nmid}
\newcommand{\abs}[1]{\lvert{#1}\rvert}
\newcommand{\partition}[1]{\ensuremath{\left\langle{#1}\right\rangle}}
\newcommand{\floor}[1]{\left\lfloor{#1}\right\rfloor}
\newcommand{\modtwo}[1]{\left\{{#1}\right\}}
\newcommand{\num}{\nu}
\newcommand{\rem}[1]{R({#1})}

\title[Races among products]{Races among products}

\author{Alexander Berkovich}
\address{Department of Mathematics, University of Florida, Gainesville,
Florida 32611-8105}
\email{alexb@ufl.edu}          

\author{Keith Grizzell}
\address{Department of Mathematics, University of Florida, Gainesville,
Florida 32611-8105}
\email{grizzell@ufl.edu}          


\date{\today}  



\begin{abstract}
In this paper we revisit a 1987 question of Rabbi Ehrenpreis. Among many things, we provide an elementary injective proof that
\[P_1(L,y,n)\geq P_2(L,y,n)\]
for any $L,n>0$ and any odd $y>1$. Here, $P_1(L,y,n)$ denotes the number of partitions of $n$ into parts congruent to $1$, $y+2$, or $2y \pmod{2y+2}$ with the largest part not exceeding $(2y+2)L-2$ and $P_2(L,y,n)$ denotes the number of partitions of $n$ into parts congruent to $2$, $y$, or $2y+1 \pmod{2y+2}$ with the largest part not exceeding $(2y+2)L-1$.
\end{abstract}

\maketitle

\section{Introduction} \label{sec:introduction}

The celebrated Rogers-Ramanujan identities \cite{RR} are given analytically as follows:
\begin{equation}
\sum^{\infty}_{n=0}\frac{q^{n^2}}{(q;q)_n} = \frac{1}{(q,q^4;q^5)_{\infty}} \label{RR1}
\end{equation}
and
\begin{equation}
\sum^{\infty}_{n=0}\frac{q^{n^2+n}}{(q;q)_n} = \frac{1}{(q^2,q^3;q^5)_{\infty}}. \label{RR2}
\end{equation}
Here we are using the following standard notations:
\begin{align*}
(a;q)_L &= \begin{cases}
              1 & \text{if } L=0,\\
              \prod_{j=0}^{L-1}(1-aq^j) &\text{if } L>0,
           \end{cases}\\
(a_1,a_2,\dots,a_n;q)_L &= (a_1;q)_L (a_2;q)_L \cdots (a_n;q)_L,\\
(a;q)_\infty &= \lim_{L\to\infty} (a;q)_L.
\end{align*}
Subtracting \eqref{RR2} from \eqref{RR1} we have
\begin{equation}
\sum^{\infty}_{n=1}\frac{q^{n^2}}{(q;q)_{n-1}}=
\frac{1}{(q,q^4;q^5)_{\infty}}- 
\frac{1}{(q^2,q^3;q^5)_{\infty}}, 
\label{RRdiff}
\end{equation}
from which it is obvious that the coefficients in the $q$-series expansion of the difference of the two products in \eqref{RRdiff} are all non-negative. In other words, for all $n>0$ we have
\begin{equation}
p_{1}(n) \geq p_{2}(n), \label{RRdiffThm}
\end{equation}
where $p_{r}(n)$ denotes the number of partitions of $n$ into parts congruent to $\pm r \pmod{5}$.

At the 1987 A.M.S.\ Institute on Theta Functions, Rabbi Ehrenpreis asked if one can prove \eqref{RRdiffThm} without resorting to the Rogers-Ramanujan identities.
In 1999, Kadell \cite{KK} provided an affirmative answer to this question by constructing an injection of partitions counted by $p_{2}(n)$ into partitions counted by $p_{1}(n)$.
In 2005, Berkovich and Garvan \cite{BG} constructed an injective proof for an infinite family of partition function inequalities related to finite products, thus giving us the following theorem.
\begin{theorem}
\label{BGThm}
Suppose $L>0$, and $1<r<m-1$. Then the
coefficients in the $q$-series expansion of the difference of the two finite products
\begin{equation*}
\frac{1}{(q,q^{m-1};q^m)_L}-
\frac{1}{(q^r,q^{m-r};q^m)_L} 
\end{equation*}
are all non-negative if and only if $r \doesnotdivide (m-r)$ and $(m-r) \doesnotdivide r$.
\end{theorem}
\noindent
We note that \eqref{RRdiffThm} is an immediate corollary of this theorem with $m=5$, $r=2$ and $L\rightarrow\infty$.

In 2011, Andrews \cite{An1} used a clever combination of injective and anti-telescoping techniques to establish the following remarkable theorem.
\begin{theorem}
\label{AnThm}
For  $L>0$, the $q$-series expansion of
\begin{equation*}
\frac{1}{(q,q^5,q^6;q^8)_L} - \frac{1}{(q^2,q^3,q^7 ;q^8)_L}
\end{equation*}
has non-negative coefficients.
\end{theorem}

The main object of the present manuscript is the following new theorem.
\begin{theorem}
\label{BerkGrizThm}
For any $L>0$ and any odd $y>1$, the $q$-series expansion of
\begin{equation}
\frac{1}{(q,q^{y+2},q^{2y};q^{2y+2})_L} - \frac{1}{(q^2,q^{y},q^{2y+1} ;q^{2y+2})_L} = \sum_{n=1}^\infty a(L,y,n) q^{n} \label{BerkGrizEqn}
\end{equation}
has non-negative coefficients. Furthermore, the coefficient $a(L,y,n)$ is $0$ if and only if either $n \in \{2,4,6,\dots, y+1\}\cup\{y\}$ or $(L,y,n)=(1,3,9)$.
\end{theorem}
We note that the products on the left of \eqref{BerkGrizEqn} can be interpreted as 
\begin{align}
\frac{1}{(q,q^{y+2},q^{2y} ;q^{2y+2})_L} &= \sum_{n=1}^\infty P_{1}(L,y,n) q^{n} \label{firstproduct}\\
\intertext{and}
\frac{1}{(q^2,q^{y},q^{2y+1} ;q^{2y+2})_L} &= \sum_{n=1}^\infty P_{2}(L,y,n) q^{n}, \label{secondproduct}
\end{align}
where  $P_1(L,y,n)$ denotes the number of partitions of $n$ into parts $\equiv 1,y+2,2y \pmod{2(y+1)}$ with the largest part not exceeding $2(y+1)L-2$ and $P_2(L,y,n)$ denotes the number of partitions of $n$ into parts $\equiv 2,y,2y+1 \pmod{2(y+1)}$ with the largest part not exceeding $2(y+1)L-1$.

In the next section, we examine a norm-preserving injection of partitions counted by $P_2(L,y,n)$ into partitions counted by $P_1(L,y,n)$, where the norm of a partition $\pi$ --- denoted $|\pi|$ --- is the sum of its parts. In section \ref{sec:proof}, we give a proof of Theorem \ref{BerkGrizThm}. In section \ref{sec:generalization}, we look at a generalization of Theorem \ref{BerkGrizThm}. In section \ref{sec:conclusion}, we conclude with a brief discussion of our plans for future work in this area.

\section{The Injection} \label{sec:injection}

Let $s_t$ denote $s + (t-1)(2y+2)$ for any positive integers $s$ and $t$, and let $\num(s_t)$ denote the number of occurrences of $s_t$ in a given partition. Then we may write any partition $\pi_1$ counted by $P_1(L,y,n)$ as
\[\pi_1 = \partition{1_1^{\num(1_1)},(y+2)_1^{\num((y+2)_1)},(2y)_1^{\num((2y)_1)},\dots,1_L^{\num(1_L)},(y+2)_L^{\num((y+2)_L)},(2y)_L^{\num((2y)_L)}},\]
where
\begin{equation}
|\pi_1| = \sum_{k=1}^L \left( 1_k\cdot\num(1_k) + (y+2)_k\cdot\num((y+2)_k) + (2y)_k\cdot\num((2y)_k) \right) = n; \label{n_for_PI1}
\end{equation}
and similarly any partition $\pi_2$ counted by $P_2(L,y,n)$ may be written as
\[\pi_2 = \partition{2_1^{\num(2_1)},y_1^{\num(y_1)},(2y+1)_1^{\num((2y+1)_1)},\dots,2_L^{\num(2_L)},y_L^{\num(y_L)},(2y+1)_L^{\num((2y+1)_L)}},\]
where
\begin{equation}
|\pi_2| = \sum_{k=1}^L \left( 2_k\cdot\num(2_k) + y_k\cdot\num(y_k) + (2y+1)_k\cdot\num((2y+1)_k) \right) = n. \label{n_for_PI2}
\end{equation}
Here we are following the convention as in \cite{An2} whereby the exponents represent the frequencies of the parts.

Let $Q(s_t)$ and $R(s_t)$ denote the quotient and remainder, respectively, upon dividing $\num(s_t)$ by $2$, taking $0 \leq R(s_t) \leq 1$.
Our injection then maps a partition $\pi_2$ to a partition $\pi_1$ as follows:
\begin{equation}
\num((2y)_k) = \begin{cases}                                              \label{BeginInjectionRules}
                Q((2y+1)_{k/2}) &\text{if $k$ is even,}\\
                Q(y_{(k+1)/2}) &\text{if $k$ is odd,}
                \end{cases}
\end{equation}
\begin{equation}
\num((y+2)_k) = \begin{cases}						\label{SecondInjectionRule}
                 \num((2y+1)_k) &\text{if $L/2 < k \leq L$,}\\
                 2\num(2_{2k})+R((2y+1)_k) &\text{if $1\leq k \leq L/2$,}
                 \end{cases}
\end{equation}
\begin{equation}
\num(1_k) = \begin{cases}
             \num(y_k) &\text{if $(L+1)/2 < k \leq L$,}\\
             2\num(2_{2k-1})+R(y_k) &\text{if $1 < k \leq (L+1)/2$,}\\
             R(y_1)+2\num(2_1)+(y-1)(A+B+C+D) &\text{if $k=1$,}
             \end{cases}
\end{equation}
where
\begin{align}
A &= \sum_{1 \leq b \leq L/2} R((2y+1)_b),                             \label{A_defined}\\
B &= \sum_{L/2 < b \leq L} \num((2y+1)_b),\\
C &= \sum_{1 \leq b \leq (L+1)/2} R(y_b),\\
D &= \sum_{(L+1)/2 < b \leq L} \num(y_b).                                  \label{EndInjectionRules}
\end{align}
For example, we have the following mappings $\pi_2 \mapsto \pi_1$ as part of our injection:
\begin{align*}
\partition{(2y+1)_k^{2i}} &\mapsto \partition{(2y)_{2k}^i} &&\text{if }1\leq k \leq L/2,\\
\partition{(2y+1)_k^{2i+1}} &\mapsto \partition{1_1^{y-1},(y+2)_k,(2y)_{2k}^i} &&\text{if }1\leq k \leq L/2,\\
\partition{(2y+1)_k} &\mapsto \partition{1_1^{y-1},(y+2)_k} &&\text{if }L/2 < k \leq L,\\
\partition{y_k^{2i}} &\mapsto \partition{(2y)_{2k-1}^{i}} &&\text{if }1\leq k \leq (L+1)/2,\\
\partition{y_k^{2i+1}} &\mapsto \partition{1_1^{y-1},1_k,(2y)_{2k-1}^i} &&\text{if }2\leq k \leq (L+1)/2,\\
\partition{y_k^{2i+1}} &\mapsto \partition{1_1^y,(2y)_1^i} &&\text{if }k=1,\\
\partition{y_k} &\mapsto \partition{1_1^{y-1},1_k} &&\text{if }(L+1)/2 < k \leq L,\\
\partition{2_k} &\mapsto \partition{1_{(k+1)/2}^2} &&\text{if $k$ is odd},\\
\partition{2_k} &\mapsto \partition{(y+2)_{k/2}^2} &&\text{if $k$ is even}.
\end{align*}

From the rules \eqref{BeginInjectionRules}--\eqref{EndInjectionRules}, it is a relatively straightforward (though perhaps slightly tedious) matter to verify that for any partition $\pi_2$ counted by $P_2(L,y,n)$, the corresponding image partition $\pi_1$ will be one that is counted by $P_1(L,y,n)$; i.e.\ if $\pi_2$ maps to $\pi_1$, then $|\pi_1| = |\pi_2|$. To show that this mapping is injective, we give the inverse map:
\begin{equation}
\num((2y+1)_k) = \begin{cases}                                              \label{BeginInverseRules}
               2\num((2y)_{2k}) + R((y+2)_k) &\text{if $1 \leq k \leq L/2$,}\\
               \num((y+2)_k) &\text{if $L/2 < k \leq L$,}
               \end{cases}
\end{equation}
\begin{equation}
\num(y_k) = \begin{cases}
               2\num((2y)_{2k-1}) + R(1_k) &\text{if $1 \leq k \leq (L+1)/2$,}\\
               \num(1_k) &\text{if $(L+1)/2 < k \leq L$,}
               \end{cases}
\end{equation}
\begin{equation}
\num(2_k) = \frac{1}{2} \cdot \begin{cases}
                               \num((y+2)_{k/2})-R((y+2)_{k/2}) &\text{if $k$ is even,}\\
                               \num(1_{(k+1)/2})-R(1_{(k+1)/2}) &\text{if $k>1$, odd,}\\
                               \num(1_1)-R(1_1)-(y-1)(W+X+Y+Z) &\text{if $k=1$,}           \label{mu:intro}
                               \end{cases}
\end{equation}
where
\begin{align}
W &= \sum_{1 \leq b \leq L/2} R((y+2)_b),\\
X &= \sum_{L/2 < b \leq L} \num((y+2)_b),\\
Y &= \sum_{1 \leq b \leq (L+1)/2} R(1_b),\\
Z &= \sum_{(L+1)/2 < b \leq L} \num(1_b).                                  \label{EndInverseRules}
\end{align}

We note that the only possibly negative quantity exhibited in either the forward map or the inverse map is the partition statistic $\mu$, which takes a partition $\pi_1$ counted by $P_1(L,y,n)$ and maps it to
\begin{equation}
\mu(\pi_1) = \num(1_1)-R(1_1)-(y-1)(W+X+Y+Z),	\label{mu:def}
\end{equation}
i.e.\ the numerator of the expression given for $\num(2_1)$ in \eqref{mu:intro}. This is a useful statistic since we must have $\mu(\pi_1) \geq 0$ if{f} $\pi_1$, counted by $P_1(L,y,n)$, is the image of some $\pi_2$ counted by $P_2(L,y,n)$. (Note that the value of $\mu$ is automatically even.)

For example, if $L=2$, $y=3$, and $n=14$, then we have the following partitions $\pi_1$ counted by $P_1(2,3,14)$ and $\pi_2$ counted by $P_2(2,3,14)$, where $\pi_2 \mapsto \pi_1$ if they are on the same row, as well as the corresponding value of the statistic $\mu$.
\begin{center}
\begin{tabular}{|ccc|}\hline
$\pi_2$ & $\pi_1$ & $\mu(\pi_1)$\\\hline
\partition{3_1,3_2} & \partition{1_1^5,1_2} & 0\\
\partition{2_1^2,2_2} & \partition{1_1^4,5_1^2} & 4\\
\partition{7_1^2} & \partition{6_2} & 0\\
\partition{2_1^2,3_1,7_1} & \partition{1_1^9,5_1} & 4\\
\partition{2_1,3_1^4} & \partition{1_1^2,6_1^2} & 2\\
\partition{2_1^4,3_1^2} & \partition{1_1^8,6_1} & 8\\
\partition{2_1^7} & \partition{1_1^{14}} & 14\\
& \partition{1_1,5_2} & $-4$\\
& \partition{5_1,1_2} & $-4$\\
& \partition{1_1^3,5_1,6_1} & $-2$\\\hline
\end{tabular}
\end{center}

\section{Proof of Theorem \ref{BerkGrizThm}} \label{sec:proof}

First, we note that the injection given in section \ref{sec:injection} proves the first part of Theorem \ref{BerkGrizThm} by virtue of the partition interpretations given by \eqref{firstproduct} and \eqref{secondproduct}. What remains to be shown is the last statement in the theorem: \emph{the coefficient $a(L,y,n)$ is $0$ if and only if either $n \in \{2,4,6,\dots, y+1\}\cup\{y\}$ or $(L,y,n)=(1,3,9)$.}

If $n \in \{2,4,6,\dots, y+1\}$, then $P_1(L,y,n)=P_2(L,y,n)=1$: here $P_1$ counts $\partition{1^n}$ and $P_2$ counts $\partition{2^{n/2}}$. If $n=y$, then $P_1(L,y,n)=P_2(L,y,n)=1$: here $P_1$ counts $\partition{1^y}$ and $P_2$ counts $\partition{y}$. If $(L,y,n)=(1,3,9)$, then $P_1(1,3,9)=P_2(1,3,9)=3$: here $P_1$ counts $\partition{1^9}$, $\partition{1^4,5}$, and $\partition{1^3,6}$; and $P_2$ counts $\partition{2^3,3}$, $\partition{2,7}$, and $\partition{3^3}$. Thus, if either $n \in \{2,4,6,\dots, y+1\}\cup\{y\}$ or $(L,y,n)=(1,3,9)$, we have $a(L,y,n)=0$.

To show the reverse implication we will show that the inverse is true: if  $n \not\in \{2,4,6,\dots, y+1\}\cup\{y\}$ and $(L,y,n)\neq(1,3,9)$, we have $a(L,y,n)>0$. To accomplish this, we use the partition statistic $\mu$ defined previously and we exhibit partitions $\pi_1$ counted by $P_1(L,y,n)$ that have $\mu(\pi_1)<0$, and hence are not mapped to under the injection. We now consider the following three cases (and several subcases), where throughout we assume that $y$ is odd.

\emph{Case 1: $y>3$ and $L\geq 1$.} Here we will examine six subcases.

\emph{Subcase 1a.}
If $0<b<2y-1$, with $b$ odd, then $\partition{1^b,(y+2),(2y)^m}$ cannot be in the image since we would have $\mu = b - 1 - (y-1)(1+1) = b - (2y-1)<0$. This implies that the inequality is strict for all even $n > y+1$ except possibly when $n \equiv y+1 \pmod{2y}$.

\emph{Subcase 1b.}
If $L=1$ then $\partition{1^{y-3},(y+2)^2,(2y)^m}$ cannot be in the image since we would have $\mu = y-3 -(y-1)(2) = -y-1 < 0$. This implies that the inequality is strict for all even $n > y+1$ with $n \equiv y+1 \pmod{2y}$ when $L=1$.

\emph{Subcase 1c.}
If $L>1$ then $\partition{1^{y-2},1_2,(2y)^m}$ cannot be in the image since we would have $\mu = y-2 - 1 - (y-1)(2) = -y-1 < 0$. This implies that the inequality is strict for all even $n > y+1$ with $n \equiv y+1 \pmod{2y}$ when $L>1$.

Note that Subcases 1a, 1b, and 1c together show that the inequality is strict for any even $n > y+1$.

\emph{Subcase 1d.}
If $0 \leq b < y-1$, with $b$ even, then $\partition{1^b,(y+2),(2y)^m}$ cannot be in the image since we would have $\mu = b - (y-1)(1) = b - (y-1)<0$. This implies that the inequality is strict for all odd $n > y$ with $n \equiv y+2,\ y+4,\ y+6, \dots,\text{ or } 2y-1 \pmod{2y}$.

\emph{Subcase 1e.}
If $0 < b < y$, with $b$ odd, then $\partition{1^b,(2y)^m}$ cannot be in the image since we would have $\mu = b - 1 - (y-1)(1) = b - y < 0$. This implies that the inequality is strict for all odd $n > 0$ with $n \equiv 1,\ 3,\ 5, \dots,\text{ or } y-2 \pmod{2y}$.

\emph{Subcase 1f.}
The partition $\partition{1^{y-4},(y+2)^2,(2y)^m}$ cannot be in the image since we would have $\mu = (y-4) - 1 -(y-1)(1) = -4 < 0$ if $L>1$ and $\mu = (y-4) - 1 -(y-1)(1+2) = -2y-2 < 0$ if $L=1$. This implies that the inequality is strict for all odd $n > y$ with $n \equiv y \pmod{2y}$.

From Subcases 1a--1f we may conclude that the inequality is strict when $y>3$, $L\geq1$, and $n\not\in\{2,4,\dots,y+1\}\cup\{y\}$.

\emph{Case 2: $y=3$ and $L>1$.}
In this case the partition $\partition{1_2,6^m}$ cannot be in the image since we would have $\mu = -2 < 0$. This implies that the inequality is strict for all odd $n > 3$ with $n \equiv 3 \pmod{6}$. Together with Subcases 1a, 1c, 1d, and 1e (all with $y=3$), this shows that if $y=3$ and $L>1$, then the inequality is strict when $n \not\in\{2,3,4\}$.

\emph{Case 3: $y=3$ and $L=1$.}
In this case the partition $\partition{5^3,6^m}$ cannot be in the image since we would have $\mu = -6 < 0$. This implies that the inequality is strict for all odd $n > 9$ with $n \equiv 3 \pmod{6}$. Together with Subcases 1a, 1b, 1d, and 1e (all with $y=3$), this shows that if $y=3$ and $L=1$, then the inequality is strict when $n \not\in\{2,3,4,9\}$.

Cases 1, 2, and 3 together show that the inequality is strict except for the following possibilities:
\begin{itemize}
\item $y>3$, $L\geq 1$, and $n \in\{2,4,\dots,y+1\}\cup\{y\}$;
\item $y=3$, $L>1$, and $n \in \{2,3,4\}$;
\item $y=3$, $L=1$, and $n \in \{2,3,4,9\}$.
\end{itemize}
However, we have already shown that the inequality is an equality at these points; thus the theorem is proven.

\section{A Further Generalization} \label{sec:generalization}

In Theorem \ref{BerkGrizThm} we may replace $y+2$ with any integer $x$, provided $1 < x \leq y+2$, and still have a perfectly viable inequality; thus, the following generalization.
\begin{theorem}
\label{BerkGrizThm2}
For any $L>0$, any odd $y>1$, and any $x$ with $1 < x \leq y+2$, the $q$-series expansion of
\begin{equation}
\frac{1}{(q,q^{x},q^{2y};q^{2y+2})_L} - \frac{1}{(q^2,q^{y},q^{2y+1} ;q^{2y+2})_L} = \sum_{n=1}^\infty a(L,y,n,x) q^{n} \label{BerkGrizEqn2}
\end{equation}
has non-negative coefficients.
Furthermore, the coefficient $a(L,y,n,x)$ is $0$ if and only if one of the following three conditions holds:
\begin{enumerate}
\item $n<x$ and $n$ is even. \label{n<x}
\item $n=y$ and $y<x$. \label{n=y}
\item $n=9$ and $(L,y,n,x)=(1,3,9,5)$. \label{n=9}
\end{enumerate}
\end{theorem}

We note that the products on the left of \eqref{BerkGrizEqn2} can be interpreted as
\begin{align}
\frac{1}{(q,q^{x},q^{2y} ;q^{2y+2})_L} &= \sum_{n=1}^\infty P_{1}'(L,y,n,x) q^{n} \label{firstproduct2}\\
\intertext{and}
\frac{1}{(q^2,q^{y},q^{2y+1} ;q^{2y+2})_L} &= \sum_{n=1}^\infty P_{2}(L,y,n) q^{n}, \label{secondproduct2}
\end{align}
where  $P_1'(L,y,n,x)$ denotes the number of partitions of $n$ into parts $\equiv 1,x,2y \pmod{2(y+1)}$ with the largest part not exceeding $2(y+1)L-2$ and, as before, $P_2(L,y,n)$ denotes the number of partitions of $n$ into parts $\equiv 2,y,2y+1 \pmod{2(y+1)}$ with the largest part not exceeding $2(y+1)L-1$.

\begin{proof}
We will prove the first part of the theorem by producing a norm-preserving injection from partitions $\pi_1$ counted by $P_1(L,y,n)$ to partitions $\pi_1'$ counted by $P_1'(L,y,n,x)$ and then relying on the fact that composition of injections is injective. Using $\num'$ to distinguish counting parts of $\pi_1'$ from counting parts of $\pi_1$, we take
\begin{equation}
\num'((2y)_k) = \num((2y)_k)
\end{equation}
\begin{equation}
\num'(x_k) = \num((y+2)_k)
\end{equation}
\begin{equation}
\num'(1_k) = \begin{cases}
              \num(1_k) & \text{if $k>1$}\\
              \num(1_1) + (y+2-x)\displaystyle\sum_{1 \leq b \leq L} \num((y+2)_b) & \text{if $k=1$}.
              \end{cases}
\end{equation}
The inverse map is immediate:
\begin{equation}
\num((2y)_k) = \num'((2y)_k)
\end{equation}
\begin{equation}
\num(x_k) = \num'((y+2)_k)
\end{equation}
\begin{equation}
\num(1_k) = \begin{cases}
              \num'(1_k) & \text{if $k>1$}\\
              \num'(1_1) - (y+2-x)\displaystyle\sum_{1 \leq b \leq L} \num'(x_b) & \text{if $k=1$}.
              \end{cases}
\end{equation}
It is then very straightforward to show that this injection is norm-preserving. Thus, we have
\begin{equation}
P_1'(L,y,n,x) - P_1(L,y,n) \geq 0,\label{P_1'>=P_1}
\end{equation}
and when we compose the injection given by \eqref{BeginInjectionRules}--\eqref{EndInjectionRules} with the one presented above, we obtain a mapping of partitions
\begin{equation}
\pi_2 \mapsto \pi_1 \mapsto \pi_1'
\end{equation}
which is an injection that maps $\pi_2 \mapsto \pi_1'$. Thus, we have
\begin{equation}
P_1'(L,y,n,x) \geq P_1(L,y,n) \geq P_2(L,y,n).
\end{equation}
For the second part of the theorem, we note that
it is straightforward to verify that for any $n$ prescribed by conditions \eqref{n<x}--\eqref{n=9}, one does, in fact, obtain $a(L,y,n,x)=0$. Also, if $P_1(L,y,n) > P_2(L,y,n)$, then $P_1'(L,y,n,x) > P_2(L,y,n)$.
So if $n$ is even and $x \leq n < y+2$, then $P_1'(L,y,n,x)\geq2$, counting at least $\partition{1^n}$ and $\partition{1^{n-x},x}$, whereas $P_2(L,y,n)=1$, counting only $\partition{2^{n/2}}$; hence condition \eqref{n<x}.
Now if $n=y$ and $x \leq y$ then $P_1'(L,y,n,x)\geq2$, counting at least $\partition{1^n}$ and $\partition{1^{n-x},x}$, whereas $P_2(L,y,n)=1$, counting only $\partition{y}$; hence condition \eqref{n=y}.
Finally, $(L,y,n,x)=(1,3,9,5)$ is the same as $(L,y,n)=(1,3,9)$ in Theorem \ref{BerkGrizThm}; hence condition \eqref{n=9}.
\end{proof}

In 1971, Andrews \cite{An3} used a simple inductive technique to prove the following theorem.
\begin{theorem}
\label{An1971Thm3}
Let $S=\{a_i\}_{i=1}^\infty$ and $T=\{b_i\}_{i=1}^\infty$ be two strictly increasing sequences of positive integers such that $b_1=1$ and $a_i \geq b_i$ for all $i$. Let $\rho(S;n)$ (resp.\ $\rho(T;n)$) denote the number of partitions of $n$ into parts taken from $S$ (resp.\ $T$). Then \[\rho(T;n) \geq \rho(S;n)\] for all $n$.
\end{theorem}
\noindent
We note that this theorem provides an alternate proof of the partition inequality \eqref{P_1'>=P_1}, as well as the subset of cases $2 \leq x \leq y$ in Theorem \ref{BerkGrizThm2}. Observe, however, that the cases when $x=y+2$ and when $x=y+1$ in Theorem \ref{BerkGrizThm2} are not covered by Andrews' theorem but are covered by our new Theorems \ref{BerkGrizThm} and \ref{BerkGrizThm2}. In addition, Theorems \ref{BerkGrizThm} and \ref{BerkGrizThm2} also provide explicit conditions for when the inequality is strict.

\section{Concluding Remarks} \label{sec:conclusion}

We plan to study more general partition inequalities in a later paper, including cases with higher modulus, cases where $y$ is even, and cases with more than three residues. Experimental evidence leads us to the following conjecture for three residues, which we are actively pursuing.
\begin{conjecture}
\label{BerkGrizConj1}
For any $L>0$, any $z>1$, any $y>z$, any $x$ with $y < x \leq y+z$, and any $m\geq yz+2$, the $q$-series expansion of
\begin{equation}
\frac{1}{(q,q^{x},q^{yz};q^{m})_L} - \frac{1}{(q^z,q^{y},q^{yz+1} ;q^{m})_L} = \sum_{n=1}^\infty a(L,y,n,x,z,m) q^{n}
\end{equation}
has only non-negative coefficients if{f} $z$ does not divide $y$; if $z$ divides $y$ then there are finitely many negative coefficients.
\end{conjecture}

In particular, we plan to prove the following.
\begin{proposal}
\label{BerkGrizProp1}
For any $L>0$ and any even $y>2$, the $q$-series expansion of
\begin{equation}
\frac{1}{(q,q^{y+2},q^{2y};q^{2y+2})_L} - \frac{1}{(q^2,q^{y},q^{2y+1} ;q^{2y+2})_L} = \sum_{n=1}^\infty a(L,y,n) q^{n}
\end{equation}
has non-negative coefficients except for $a(L,y,y) = -1$.
\end{proposal}
Note that Conjecture \ref{BerkGrizConj1} with $z=2$, $m=2y+2$, and $y$ odd is part of Theorem \ref{BerkGrizThm2}, and that Proposal \ref{BerkGrizProp1} is the natural companion to Theorem \ref{BerkGrizThm}. Also note that, as before, Theorem \ref{An1971Thm3} would clearly establish the corresponding result to Conjecture \ref{BerkGrizConj1} when $x\leq y$, leaving the more difficult cases when $x>y$ to be addressed.

Finally, we would like to point out that the problems discussed in this paper belong to a broad class of positivity problems in $q$-series and partitions. These problems often are very deceptive because they are so easy to state but so painfully hard to solve. As an example, consider the famous Borwein problem:
\begin{quote}
Let $B_e(L,n)$ (resp.\ $B_o(L,n)$)  denote the number of partitions of $n$ into an even (resp.\ odd) number  of distinct nonmultiples of $3$ with each part less than $3L$. Prove that for all positive integers $L$ and $n$,
$B_e(L,n)-B_o(L,n)$  is nonnegative if $n$ is a multiple of $3$ and nonpositive otherwise.
\end{quote}
Further background on this conjecture may be found in \cite{An4}, \cite{BW}, \cite{BR}, \cite{W1}, and \cite{W2}.

\vspace{1em}

\noindent
{\it Acknowledgements}.
We are grateful to George Andrews for bringing our attention to \cite{An3}, and to Krishna Alladi and Frank Garvan for their interest.

\providecommand{\bysame}{\leavevmode\hbox to3em{\hrulefill}\thinspace}

\end{document}